\documentstyle{article}
\title{Vanishing via lifting to second Witt vectors  \\and
 a   proof of an isotriviality result}
\author{Mark Andrea A.  de Cataldo\thanks{
Partially supported by N.S.F. Grant DMS 9701779.}}

\date{August 1998}

\newtheorem{tm}{Theorem}[section]
\newtheorem{lm}[tm]{Lemma}
\newtheorem{pr}[tm]{Proposition}
\newtheorem{rmk}[tm]{Remark}
\newtheorem{cor}[tm]{Corollary}

\newtheorem{??}[tm]{Question}

\font\tenmsb=msbm10
\font\sevenmsb=msbm7
\font\fivemsb=msbm5

\newfam\msbfam
\textfont\msbfam=\tenmsb
\scriptfont\msbfam=\sevenmsb
\scriptscriptfont\msbfam=\fivemsb
\def\Bbb#1{{\fam\msbfam #1}}

\font\teneufm=eufm10
\font\seveneufm=eufm7
\font\fiveeufm=eufm5
\newfam\eufmfam
\textfont\eufmfam=\teneufm
\scriptfont\eufmfam=\seveneufm
\scriptscriptfont\eufmfam=\fiveeufm
\def\frak#1{{\fam\eufmfam\relax#1}}

\newcommand\ci{\cite}
\newcommand\s{\sigma}

\newcommand\zed{{\Bbb Z}}

\newcommand\pn[1]{{\Bbb P}^{#1}}

\newcommand\blacksquare{{\hspace*{\fill} $\Box$}} 
\newcommand\odix[1]{ {\cal O}_{#1} }

\newcommand\surj{-\!\!\!-\!\!\!-\!\!\!-\!\!\!\!\!\gg}

\newcommand\om[2]{\omega_{{#1}/{#2}}}
\newcommand\ff[1]{\frak {#1}}
\newcommand\omk[3]{\omega^{\otimes {#3}}_{{#1}/{#2}}}

\begin{document}
\maketitle
\begin{abstract}
A  proof based on reduction to finite fields
 of Esnault-Viehweg's 
stronger version of  Sommese Vanishing Theorem for $k$-ample line  bundles
is given.
This result is used to  give different   proofs of   isotriviality
results of A. Parshin and L. Migliorini. 
\end{abstract}

\section{Introduction}
\label{intr}
This note contains a   proof of  Esnault-Viehweg's improvement of 
 Sommese Vanishing Theorem for $k$-ample line bundles.
The proof is based on   M. Raynaud's
proof of  Akizuki-Kodaira-Nakano Vanishing Theorem.

A vector bundle version and a Weak-Lefschetz-type Theorem, which are easy
consequences of the vanishing result, but do 
not seem to have appeared in the literature, are also given.
 
The vanishing theorem allows to give an algebraic proof
of a result of A. Parshin and L. Migliorini on the isotriviality 
of smooth fibrations over curves of genus at most one
with fibers either curves of genus at least two, or minimal surfaces of general type. 

\medskip
The paper is organized as follows.
\S\ref{s-o} contains, for the convenience of the reader, basic known facts about the spreading out technique. 
\S\ref{vancharp} contains the proof of the vanishing in the line bundle case, 
Theorem \ref{mine} and of its easy corollaries  Corollary \ref{vb} and Corollary \ref{weaklef}. 
\S\ref{pm} contains a proof of
A. Parshin and L. Migliorini's result:
Theorem \ref{tpm}.

\medskip
\noindent
{\bf Acknowledgment.}
This paper was written while the author enjoyed the hospitality of the Max-Planck-Institut
f\"ur Mathematick in Bonn.

\section{Preliminaries on the ``spreading out" technique.}
\label{s-o}
Any reasonable and finite amount of geometry defined over a field
of characteristic zero $ K $
can be ``spread out" over an algebra of finite type  over $\zed$, 
$ R$,  contained in $ K $.
 We can then try to use the resulting ``fibration" to say something about
the original situation over $K$.

This very vaguely presented principle is made precise in EGA, IV 8, 10.
A concise exposition  is given by L. Illusie in \ci{il},  \S6.

\medskip
Let us exemplify the ``spreading out" procedure by listing the properties we need
in the sequel of the paper.

\smallskip
Let us fix some notation.
 Let  $Z$ be a scheme, $g: U\to V$ be a $Z$-morphism of $Z$-schemes, $F$ be
a $\odix{U}$-module and $z$ be a point in $Z$. We denote the fibers 
$U \times_Z Spec\, k(z)$ and $V \times_Z Spec \, k(z)$ simply by
$U_z$ and $V_z$, the restriction of $F$ to $U_z$ by $F_z$, and the induced morphism
from $U_z\to V_z$ by $g_z$. 

\noindent
A line bundle is said to be {\em semi-ample} if some positive power of it is generated by its global sections. Note that this does not
imply  the (false) statement that every sufficiently high power of it is generated by its global sections; e.g. a non-trivial torsion
line bundle over an elliptic curve.

\noindent
Let $b$ be a non-negative integer; a
 semi-ample line bundle $L$ on  a projective variety
$X$ is said to be {\em $b$-ample} if, given
a positive integer $N$ such that $L^{\otimes N}$ defines a morphism
$\phi_{|NL|}: X \to {\Bbb P}$, then the non-empty fibers of $\phi_{|NL|}$
have dimension at most $b$. This notion does not depend
on $N$. Note that the Kodaira-Iitaka dimension
$\kappa (L)= \dim{\phi_{|NL|} (X)}$ for any $N$ as above;
see \ci{mo2}.

\noindent
Given a ring $T$, $W_2(T)$ is the ring of Witt vectors
of length two associated with $T$.
\begin{pr}
\label{spread}
Let $f:X \to Y$ be a projective $ K $-morphism of projective $ K $-varieties,
$F$, $L$ and $A$, respectively, be a coherent sheaf on $X$, a $b$-ample line bundle
on $X$
and an ample line bundle on $Y$, respectively.

\smallskip
\noindent
Then, there exist an integral $\zed$-algebra of finite type, $R$, contained in 
$ K$,
projective 
$Spec \,  R= : {\cal R}$-schemes $\rho : {\cal X} \to {\cal R}$ and $\s: 
{\cal Y} \to {\cal R}$, a projective $\cal R$-morphism 
$\ff{f}: {\cal X} \to {\cal Y}$, coherent sheaves ${\cal F}$ and ${\cal L}$
on ${\cal X}$ and ${\cal A}$ on ${\cal Y}$, and a Zariski-dense open
subset  ${\cal U} \subset\cal R$, contained
in the  locus of $\cal R$ which is smooth over $Spec\,\zed$, with the properties listed below.

(i)
The objects $X$, $Y$, $f$, $F$, $L$ and $A$, respectively, are obtained
from the corresponding objects ${\cal X}$, ${\cal Y}$, $\ff{f}$, ${\cal F}$,
${\cal L}$ and ${\cal A}$, respectively, by means of the base change induced by
$ R \hookrightarrow  K $.

(ii)
The objects $\rho$, $\s$, ${\cal F}$, ${\cal L}$, ${\cal A}$ and $R^i \ff{f}_* {\cal F}$
are all flat over $\cal U$. In particular the formation of the sheaves
$R^i \ff{f}_* {\cal F}$ commutes with taking the fiber over any point
$u \in {\cal U}$.

(iii) The sheaves ${\cal L}$ and ${\cal A}$ are locally free on
$\rho^{-1} ({\cal U})$ and $\s^{-1} ({\cal U})$, respectively,
${\cal A}$ is  $\rho$-ample and ${\cal L}_u$ is $b$-ample for every
$u \in {\cal U}$.

(iv) If $X$ is smooth,  then $\cal U$ can be chosen so that $\rho$  and $\s$,
respectively, are smooth and flat, respectively, over ${\cal U}$. 

(v) If $s$ is a closed point in $\cal U$, then
${\cal X}_s$ lifts to $W_2(k(s))$. Moreover, we can choose
$\cal U$ so that
$char \, k(s) > \dim {X}$, for every closed point $s \in {\cal U}$.
\end{pr} 
{\em Proof.} The first set of properties follows from
EGA IV 8. The second one is the theorem on flattening stratifications 
of Grothendieck, and cohomology and base change as established in
EGA III 6.9.10. The third is EGA IV 8; the statement  about 
$b$-ampleness stems from the fact that
a line bundle $L$ on $X$ is $b$-ample and not
$(b-1)$-ample iff there exists  ample divisors
$D_1, \ldots, D_b$ such that $L_{|\cap_{i=1}^b D_i}$
is ample and $b$ is the minimum number for which this is possible.
 The fourth one follows from  generic smoothness and generic flatness, respectively.
The last one follows  after  shrinking $\cal U$, if necessary,
so that the conclusion on the characteristics is true,
and  by the fact that $\cal U$ is smooth (over $\zed$); 
see \ci{il}, page 152-3. 
\blacksquare

\medskip
The following elementary lemma contains basic facts to be used later.
\begin{lm}
\label{densw2}
Every closed point in  ${\cal R}$ has finite, and a fortiori  perfect, residue field.

\noindent
Every Zariski-dense open subset of $\cal R$ contains a Zariski-dense set
of closed points.

\end{lm}
{\em Proof.} EGA IV 10.4.6, 10.4.7
\blacksquare

\medskip
 The following   result   contains the basic fact  that we shall need
 about smooth projective varieties
 over a perfect field $\frak k$ which lift to $W_2({\frak k})$.
It is proved in \ci{de-il} as a consequence of Th\'eor\`eme 2.1 (and Corollaire 2.3).  
\begin{tm}
\label{dikan}
(Akizuki-Kodaira-Nakano Vanishing Theorem) Let $X$ be a smooth projective variety of dimension $d$
 over a perfect field  $\frak k$ of characteristic
$p> d$ which admits a lifting to $W_2({\frak k})$, 
$M$ be a line bundle on $X$ and
 $\nu $ be a non-negative  integer.

\noindent
Assume that, for some positive integer $n$, 
$$
H^j(X, \Omega_X^i \otimes M^{\otimes p^n})=\{0 \} \qquad \forall \, i+j =\nu.
$$
Then
$$
H^j(X, \Omega_X^i \otimes M)=\{0 \} \qquad \forall \, i+j =\nu.
$$ 
In particular, if $\check{M}$ is ample, then any $\nu < d$  will do. 
\end{tm}
{\em Proof.} See \ci{de-il} Corollaire 2.3 and  Lemme 2.9. Note that
in our setting we can conclude that the relevant hypercohomology group
${\Bbb H}^{\nu}=\{0 \}$, which is what is needed. \blacksquare

\medskip
If $M$ is ample, then one proves Kodaira-Akizuki-Nakano Vanishing Theorem
as an easy consequence
of Serre Vanishing (Raynaud).
This theorem allows one to re-prove the classical
Kodaira-Akizuki-Nakano Vanishing Theorem  in characteristic zero
by spreading out $X \to  K $ to ${\cal X} \to {\cal R}$, $\Omega^{\bullet}_{X/ K }$
to $\Omega^{\bullet}_{ {\cal X}/{\cal R} }$, by the upper semi-continuity properties of the dimensions of cohomology groups.

\section{Esnault-Viehweg's improvement of A. Sommese Vanishing Theorem
for $k$-ample vector bundles  using $W_2$-lifting.
}
\label{vancharp}
The following result is a slight improvement of A. Sommese Vanishing Theorem
for $k$-ample line bundles; see \ci{s-s}, Theorem 3.36 and Corollary 5.20. This improvement in the line bundle case is due to 
H. Esnault and E. Viehweg, \ci{e-v1}, Theorem 2.4, who use analytical
methods. 

\smallskip
We offer a new proof which is algebraic and passes through reduction to finite fields, 
Deligne-Illusie decomposition Theorem and  characteristic zero
vanishing theorems made valid in the finite field case by ``propagation." 
 This technique has been already employed in the context of non-complete
varieties in \ci{b-k} where one finds,  as a particular case, a proof of
 Sommese Vanishing Theorem.

\smallskip
The ``shifted" version for the  vector bundle case follows, as it is 
nowadays standard, by 
a theorem of  J. Le-Potier's simplified by M. Schneider; see \ci{s-s}, Theorem 5.16.   

\medskip
We shall need the following fact  in the sequel when we shall
need to use the Improved Grauert-Riemenschneider Theorem
in the finite field context, where it does not hold in general.
\begin{lm}
\label{grauert}
Let $f:X\to Y$, $F$, ${\cal X} \stackrel{\ff{f}}\to {\cal Y}$ and 
${\cal F}$ be  as 
 in Lemma \ref{spread} and $\eta$ be the generic point
of $\cal R$.
Assume that
$R^j f_*  F=0$ for a fixed integer $j$. 

\noindent
Then $R^j {\ff{f}_{u}}_* {\cal F}_{u}=0$
for every point $u$ in  a suitable Zariski-dense open subset $\cal U$ of $\cal R$.
\end{lm}
{\em Proof.}
The base changes induced by $  R \hookrightarrow k( \eta )$ and $k(\eta )
\hookrightarrow  K $ are
 both flat, the second one is even faithfully flat. It follows that
 $ [R^j {\ff{f}}_* {\cal F}]_{\eta}= R^j {\ff{f}_{\eta}}_* {\cal F}_{\eta}=0$.
Since $\s$ is proper, we see that, after shrinking $\cal U$ if necessary, 
the assertion holds by Lemma \ref{spread}.ii.
 \blacksquare

\medskip
We fix $ K $,  a field of characteristic zero.
\begin{tm}
\label{mine}
Let $X$ be a smooth projective variety of dimension $d$
over $ K$, $L$ be a 
$b$-ample line bundle over $X$ of Kodaira-Iitaka dimension
$\kappa (L)$.

\noindent
Then 
\begin{equation}
\label{3}
H^j(X, \Omega^i_X \otimes L^{\vee})=\{0 \} \qquad \forall \, (i,j)
\, \, s.t.\quad 
i+j ~ < ~ \min{({\kappa (L)}, d - b +1)}. 
\end{equation}
Equivalently,
\begin{equation}
\label{4}
H^q(X, \Omega_X^p \otimes L)=\{0 \} \qquad \forall \, (p,q) \, \,  s.t.
 \quad 
p+q ~ > ~  2d - \min{({\kappa (L)}, d - b  +1)}.
\end{equation}
\end{tm}
{\em Proof.} The two statements are equivalent by virtue of  Serre Duality and of the 
canonical isomorphisms ${\Omega^{d-l}_X} \simeq 
K_X \otimes \check{\Omega}^l_X$.

\medskip
\noindent
{\em STEP I.}
We first prove the assertion under the additional hypothesis
that $mL$ is generated by its global sections for
{\em every} $m\gg 0$.

We prove  statement (\ref{4}).

There are  a surjective and projective morphism $g:X \to Y$
with connected fibers onto a normal variety $Y$
and an ample line bundle $A$ on $Y$ such that $L=g^*A$.

\noindent
By assumption, $Y$ has dimension $\kappa (L)$ and the fibers of
$g$ are at most $b$-dimensional.

\medskip
\noindent
We have the following two properties:

(a) $R^l g_{*}$ is the zero functor, for every $l > b$;

(b) $R^l g_{*} K_X =0$ for every  $l > d - \kappa (L)$; this Improved Grauert-Riemenschneider Vanishing Theorem  follows from the improved 
Kawamata-Viehweg Vanishing Theorem, \ci{s-s} Corollary 7.50  and
 from  \ci{c-k-m}, Proposition 8.9. Note that everything is algebraic here
and that H. Hironaka's Resolution of Singularities is needed.

\medskip

 We apply  Proposition 
\ref{spread}  to $g:X \to Y$, $L$, $A$ and $F:=K_X\simeq
\om{X}{ K }$ with the choice  ${\cal F}\simeq \om{ {\cal X} }{\cal R}$.
Let $s \in {\cal R}$ be a closed point  belonging to the open set
$\cal U$
(recall Lemma \ref{densw2}!) over which
all the conditions of Proposition \ref{spread}  for $\om{\cal X}{\cal R}$
and all of its direct images via $\ff{g}$
 are met. 

By virtue of Lemma \ref{grauert}, we may shrink
$\cal U$, if necessary, so that the two  conditions  (a) and (b) above are met for $\ff{g}_s$ as well.

By abuse of notation, we denote ${\cal L}_s$,  ${\cal A}_s$, ${\cal X}_s$, ${\cal
Y}_s$ and $\ff{g}_s$,
respectively, by $L_s$,  $A_s$, $X_s$. $Y_s$ and $g_s$.
In order to apply Theorem \ref{dikan}, we need to check that
$H^q(X_s, \Omega^p_{X_s} \otimes L^{\otimes m}_s)=\{0\}$, for every $m \gg 0$ in the prescribed range for $p$ and $q$. 
This is an immediate consequence
of the Leray spectral sequence for $g_s$ and conditions (a) and (b) for $g_s$,
as we now show. 

By virtue of  Leray spectral sequence and of  Serre Vanishing,
there exists $m_s$ such that for every $m\geq m_s$ and for every pair of indices
$p$ and $q$ we have that 
$$
H^q(X_s, \Omega_{X_s}^p \otimes L^{\otimes m}_s ) \simeq H^0(Y_s, R^q {g_s}_{*}
\Omega_{X_s}^p \otimes A^{\otimes m}_s).
$$
If  $p$ and $q$ are in the prescribed range, then 
$R^q {g_s}_{*} \Omega_{X_s}^p$ is zero so that the  group  on the right vanishes: in fact, $p+q \geq d +b$ and $p\leq d$ so that  we can use (a) and (b) above.

\medskip
Now the standard semi-continuity argument. 
By virtue of what we have proved and by virtue of  Proposition \ref{spread}
we can assert the existence  of  a Zariski-dense open subset 
 ${\cal U} \subset {\cal R}$ over which
the sheaves ${\cal L}$ and $\Omega^{\bullet}_{{\cal X}/{\cal R}}$ are locally free
 and  such that, given any closed point
$s \in {\cal U}$, we have that $char \, k(s) > d$,
$X_s$ lifts to $W_2(k(s))$  and  there is  a certain positive integer
$m_s$ such that for every $p$ and $q$ in the prescribed range
and for every $m\geq m_s$ we have that 
$$
H^q (X_s, \Omega_{X_s}^p \otimes L_s^{\otimes m}) = \{ 0 \}.
$$
We choose $m:=[char \,  (k(s)]^{m_s}$ and  apply  a straightforward
descending induction coupled with
Theorem  \ref{dikan} to deduce that the vanishings above hold with $m=1$.

\medskip
\noindent
The vanishing in characteristic zero follows from the fact that, by the 
upper-semiconti\-nuity
of these dimensions, we have the vanishing over the generic point
$\eta \in {\cal R}$ and, by  the flat base change
induced by $k(\eta) \hookrightarrow  K$, therefore over $ K$.

\medskip
The proof of {\em STEP I} is now complete.

\medskip
\noindent
{\em STEP II}. We now remove the additional assumption of {\em STEP I}:
we prove the theorem
by induction on $d- \kappa (L)$ using {\em STEP I} and by means of an easy
procedure  to construct, 
on a suitable covering of $X$,
 $d- \kappa (L)$ sections of the pull-back of  $L$
with base locus of dimension $d- \kappa (L)$.

Note that we may assume that $\kappa  (L) >0$, since, if $\kappa (L) =0$, 
then there is nothing to prove. 

Let $c$ be a positive  integer such that $cL$ is generated by its 
global sections. We get a surjective morphism onto a variety
$Y$ of dimension $\kappa (L)$. Choose a general section
$\s_1$
of $cL$ so that it defines a smooth divisor $D'_1$ on $X$.
Consider the corresponding cyclic covering $C_1: X_1 \to X$
branched along $D'_1$ and ramified along the smooth divisor
$D_1:= C_1^{-1}(D'_1)$ which is the zero set of  a section
of  the line bundle $L_1:=C_1^*
L$. We also have that $\Omega^l_X$ is a direct summand
of  ${C_1}_* \Omega^l_{X_1}$  for every integer $l$
such that $0 \leq l \leq d$ (see \ci{e-v}, page 6).

By iterating this procedure, we obtain a sequence of cyclic coverings
$X_{i+1} \stackrel{C_{i+1}}\longrightarrow X_{i}$ together
with line bundles $L_{i+i} = C_{i+1}^* L_{i}$ 
such that $\dim{Bs |L_{ \kappa (L)}|} \leq d- \kappa (L)$
and $\Omega^l_X$ is a direct summand 
of $C_*\Omega^l_{X_{\kappa (L)}}$, where $C:X_{\kappa (L)}
\to X$ is the induced morphism.
Since $C$ is finite, 
$H^r(X, \Omega_X^l \otimes L)$
is a direct summand of $H^r(X_{\kappa (L) },
\Omega^l_{X_{\kappa (L)}} \otimes L_{\kappa (L) } )$.

It is therefore  enough to prove the theorem
under the additional assumption that $|L| \neq \emptyset$
and that
$\dim{ Bs |L|  }\leq d- \kappa (L)$.

\smallskip
We work by  ascending induction
on ${\frak l} : =d- \kappa (L)$.

Let ${\frak l} =0$. Then $\dim {Bs |L|} =0$. By  
Zariski's \ci{za}, Theorem 6.2, we obtain that 
$mL$ is generated by its global sections for {\em every} $m\gg 0$
and we conclude by virtue of  {\em STEP I}.

Let the contention be true for every
${\frak l}' < {\frak l}$. Let us prove it for 
${\frak l}$. 
We prove   statement  (\ref{3}). 

Let $H$ be an ample hypersurface of $X$ such that
it is smooth, $L_{|H}$ is $(b-1)$-ample and $\kappa (L_{|H})=
\kappa (L) $. Such an $H$ exists by  a result of H. Hironaka's,
\ci{s-s}, Theorem 3.39,
and by the fact that $H$ maps onto $Y$ by the assumption
$d- \kappa (L) >0$. 
Consider the exact sequences:
\begin{equation}
\label{1}
0 \to \Omega^l_X \otimes L^{\vee} \otimes H^{\vee} 
\to
\Omega ^l_X \otimes L^{\vee} 
\to
\Omega ^l_X \otimes L^{\vee} 
\otimes \odix{H}  \to 0,
\end{equation}
\begin{equation}
\label{2}
 0 \to
\Omega^{l-1}_H \otimes L^{\vee}_{|H}  \otimes H^{\vee}_{|H}
\to
\Omega^l_X \otimes L^{\vee} \otimes \odix{H} 
\to
\Omega^l_H \otimes  L^{\vee}_{|H} 
\to
0.
\end{equation}
Since $L+H$ is ample, we can use the standard
Akizuki-Kodaira-Nakano Vanishing Theorem on $X$ and on $H$
and get, for $i$ and $j$ such that
$i+j < \min{(\kappa(L),d-b+1)}$, the following two 
injective maps
$$
 H^j(X, \Omega^i_X \otimes L^{\vee}  )
\hookrightarrow
 H^j(H, \Omega^i_X \otimes L^{\vee} \otimes \odix{H})
\hookrightarrow
 H^j (H, {\Omega^i_H} \otimes L^{\vee}_{|H}).
$$
Since
$$
{\frak l}_H:= ~(d-1) - \kappa (L_{|H}) ~ = ~ (d-1) - \kappa (L)  ~ < ~ d - \kappa (L) ~ = ~
{\frak l}
$$
and 
$$
i+j ~ < ~ \min{ (\kappa (L), d - b +1)} ~ = ~
\min { ( \kappa (L_{|H}) , (d-1 - (b-1) +1 ) },
$$
 we can apply the induction hypothesis and conclude that the last group on the right is trivial.
This gives the wanted vanishing result.
\blacksquare 

\begin{rmk}
\label{slight}
{\rm
Sommese's theorem is sharp, as stated. However,
note that $\kappa (L) \geq d-b$ and that the strict inequality is possible.
 If we have equality, then $f$ is equidimensional and we get Sommese's statement. If we have strict inequality,
then Theorem \ref{mine} improves Sommese's   by one unit.

\noindent
Consider the case where 
$K$ is algebraically closed and $f: X \to Y$ is the blowing-up of $\pn{3}$ at either a (closed) point, or 
along a line. Let $L$ be the pull-back of the hyperplane bundle.
 In the former case $\min{ (\kappa (L), d-b+1) }=\min{(3,2)}=2$; 
Sommese Vanishing Theorem  predicts vanishing for 
$1+j < 1$;  Theorem \ref{mine} predicts vanishing for $i+j< 2$; moreover,
$H^1(X,\Omega_X^1 \otimes  L^{\vee})$ is one dimensional.  
In the latter case $\min{ (\kappa (L), d-b+1) }=\min{(3,3)}=3$; 
A. Sommese's Theorem predicts vanishing for $i+j < 2$ and
Theorem \ref{mine} for $i+j <3$.

\noindent
This example shows that Theorem \ref{mine} is sharp and improves upon A. Sommese's. Moreover it shows concretely why
it is sharp: for in the case we blow-up a point  $p$, we have
that $R^1f_* \Omega_X^1$ is isomorphic to the skyscaper sheaf at
$p$ of stalk $K=k(p)$. The second case   puts in evidence
that, for the purpose of Akizuki-Kodaira-Nakano-type statements,
a line bundle which is 
semi-ample, big and $1$-ample is as good as ample; we use this fact
in an essential way in the proof of Theorem \ref{tpm}.
}
\end{rmk}

\medskip

The following two results follow easily from Theorem \ref{mine} 
and they do not need
the  proof given above. The first one admits a dual formulation which we  omit
for brevity. The interested reader can easily formulate and prove
vanishing results analogous to the first one which involve
$K_X \otimes \wedge^l E$ and more generally $\Omega_X^p\otimes \wedge^l E$;
see \ci{s-s}, \S5. 

\smallskip
Recall that a vector bundle $E$ is said to be $b$-ample if the associated tautological line bundle $\xi_E$ is $b$-ample.
\begin{cor}
\label{vb}
Let things be as in Theorem \ref{mine} except that we replace
the line bundle $L$  by  a rank
$r$, $b$-ample vector bundle $E$ and we set $\kappa (E) : = \kappa (\xi_E)$.

\noindent
Then 
$$
H^q (X, \Omega^p_X \otimes {E})=\{ 0 \}  \qquad \forall \, (p.q)~s.t. \quad p+q  >
 2 \, [d+(r -1)] - \min{(\kappa (E) , d+ r -b ) }
$$
\end{cor}
{\em Proof.}
By virtue of a result of Le Potier's (cf. \ci{s-s}, 5.17, 5.21 and 5.28):

\noindent
 $H^q(X, \Omega_X^p \otimes E)\simeq
H^q ( {\Bbb P}(E), \Omega_{{\Bbb P}(E)}^p \otimes \xi_E)$. Apply Theorem
\ref{mine} to the pair $( {\Bbb P} (E) ,  \xi_E )$.
 \blacksquare

\begin{cor}
\label{weaklef}
(Weak Lefschetz Theorem)
Let $X$ be a smooth projective variety of dimension
$d$ defined over $K$.
Let $D$ be an effective smooth divisor on $X$
such that the associated line bundle $L$ is $b$-ample and it is of Kodaira-Iitaka
dimension $\kappa (L)$.

\noindent
Then the canonical morphisms of de Rham cohomology
$H^l_{DR}(X/K) \to H^l_{DR} (D/K)$ are:

(i) isomorphisms for $l< \min {(\kappa (L), d-b+1)} -1$,

(ii) injective for $l < \min{ (\kappa (L), d-b +1)} $.
\end{cor}
{\em Proof.} Note that $\kappa (L_{|D}) =
\kappa (L) -1$ and that $b_D \leq  b$. We conclude by means of
easy diagram considerations on the
long cohomology sequences  associated with (\ref{1}) and (\ref{2}).
\blacksquare

\smallskip
For more statements in the vein of the corollary above, see \ci{s-s}, Theorem 3.40.

\section{A proof of  a result of A.N. Parshin and  L. Migliorini}
\label{pm}
We give an algebraic proof  of Theorem
\ref{tpm} below.
The proof hinges on Theorem \ref{mine} and on the following
positivity result of J. Koll\'ar (which holds in  greater generality
than the one stated below), whose original proof is  Hodge-theoretic
and  which has been 
proved again algebraically by
J. Koll\'ar and E. Viehweg.

\smallskip
In what follows everything is defined over   an algebraically 
closed field of characteristic zero $K$.

\begin{tm}
\label{kmt}
(Cf. \ci{ko87})
Let $f: X \to P$ be a surjective morphism with connected fibers
 of nonsingular projective varieties, where
$P$ is  a nonsingular curve. Assume that the fibers of $f$ are of general type 
and are not all birationally isomorphic to each other.

\noindent
Then the vector bundle $f_* \omk{X}{P}{m}$ is ample for infinitely many values
of the positive integer $m$. 

\end{tm}
\begin{tm} 
\label{tpm}
(Cf. \ci{pa}, \ci{mi}) Let $X$ be a
nonsingular projective variety of dimension $d$,
 $P$ be a nonsingular complete  curve of genus $g(P)\leq 1$ and
 $f:X \to P$ be a surjective smooth morphism  such that all the fibers are connected varieties
of general type with nef canonical bundle.

\noindent
If $d\leq 3$, then all the fibers are isomorphic to each other.
\end{tm}
{\em Proof.}
Since, if necessary, we  can take a double cover $P \to \pn{1}$, $P$
any elliptic curve,  we can assume that $g(P)=1$.
Note that, in this case,  $K_X = \om{X}{P}\simeq \Omega^{d-1}_{X/P}$.

Seeking a contradiction, we assume that the fibers of $f$ are not all 
birationally isomorphic to each other.

By the Base-Point-Free Theorem of Y.  Kawamata and V. Shokurov (cf. 
\ci{c-k-m}) applied to the pluricanonical line bundles of the fibers and by  Noetherian induction on $P$, there exists a positive  integer $m_0$
such that:  for every $m\geq m_0$, the natural morphism
$f^*f_* mK_X \to mK_X$ is surjective and it induces a $P$-morphism
$g_m: X \to {\Bbb P}(f_* \omk{X}{P}{m})$ with the property that
$m K_X\simeq g^*_m \xi_{E_m}$.  This morphism induces the
birationally isomorphic stable pluricanonical morphisms
on the fibers of $f:X \to P$.

By virtue 
of Theorem \ref{kmt},   we can choose the integer $m$ above so that
$E_m:= f_* \omk{X}{P}{m}$ is ample.
 It follows that $m K_X$, being the pull-back of an ample line bundle
via $g_m$, is semi-ample and $(d-2)$-ample. This conclusion holds for $K_X$ as well.

\smallskip
The  following argument is due to
  S. Kov\'acs (cf. \ci{ko}, page 370). 
Consider the exact sequences:
$$
0 \to \Omega^{i-1}_{X/P} \otimes K_X \to \Omega^i_X \otimes K_X
\to
\Omega^i_{X/P} \otimes K_X \to 0.
$$ 
For every $1\leq p \leq d-1$, we get short exact sequences
$$
H^{d-p}(X, \Omega^p_{X/P} \otimes K_X ) 
\stackrel{\alpha_p}\longrightarrow 
H^{d-(p-1)}(X, \Omega^{p-1}_{X/P} \otimes K_X)
\longrightarrow
H^{d
-(p-1)}(X,\Omega^p_X \otimes K_X),
$$
where, when $d\leq 3$,  the maps $\alpha_p$ are all surjective by Theorem \ref{mine}.
We compose all these surjective maps $\alpha_p$ and get
a   surjection
$$
\{ 0 \}= H^1(X, K_X \otimes K_X) \surj
H^d(X, K_X)\simeq { K},
$$
the first isomorphism on the left being Kawamata-Viehweg
Vanishing Theorem. This is a contradiction.

The fibers are therefore birationally isomorphic to each other. 
The result follows from the uniqueness of minimal models for curves and surfaces.
\blacksquare

\begin{rmk}
{\rm
The case $d=2$ is due to A.N. Parshin; see \ci{pa}. The case
$d=3$ is due to L. Migliorini; see \ci{mi}; his proof uses
analytic techniques.

\noindent
 The result is false without the restriction on the genus of the base; see \ci{kod}.

\noindent
Since
the automorphism group of the fibers is finite and
the fibers are all isomorphic to each other,  the fibration is isotrivial, i.e. it becomes trivial after a finite base change
$P'\to P$. 

\noindent
If we drops the nefness assumption, then
a birationally isomorphic statement still holds; see  \ci{mi}.

\noindent
A similar but weaker statement holds for any
value of $d$ and it is due to 
S. Kov\'acs \ci{ko} who has also proved the case $d=4$ of Theorem \ref{tpm}. 

\noindent
Q. Zhang \ci{qz} has proved a similar statement in any dimension
under the assumption that all fibers have ample canonical bundle.
}
\end{rmk}



\bigskip
\noindent
Author's address:
Department of Mathematics, Harvard University, One Oxford Street, 
Cambridge, MA 02138 $\quad$
e-mail: {\em mde@math.harvard.edu}


\begin{thebibliography}{99}

\bibitem{b-k} I. Bauer, S. Kosarew, ``Kodaira vanishing theorem on
 noncomplete algebraic manifolds," Math. Z. {\bf 205} (1990), no.2, 223-231.


\bibitem{c-k-m} H. Clemens, J. Koll\'ar, S. Mori, 
{\em Higher dimensional complex geometry}, Ast\'erisque {\bf 166} (1988)


\bibitem{de-il} P. Deligne, L. Illusie, ``Rel\`evetments
modulo $p^2$ et d\'ecomposition du complexe de de Rham," Invent.
math. {\bf 89}  (1987), 247-270.



\bibitem{e-v1} H. Esnault, E. Viehweg, ``Vanishing and nonvanishing theorems,"
Ast\'erisque {\bf 179-180} (1989), 97-112.



\bibitem{e-v} H. Esnault, E. Viehweg, ``Lectures on vanishing theorems,"
DMV Seminar, Band 20, Birkh\"auser, Boston (1992)



\bibitem{il} L. Illusie, ``Frobenius et d\'eg\'en\'erescence de Hodge," in J. Bertin, J.-P. Demailly, L. Illusie,
C. Peters, {\em Introduction \`a la Th\'eorie de Hodge}, Panoramas
et Synth\`eses, N{\bf 3}, Soc. Math. de France  (1996).



 


\bibitem{kod} K. Kodaira, ``A certain type of irregular surface,"
Jour. d' Analyse Math. {\bf 19} (1967).


\bibitem{ko87} J. Koll\'ar, ``Subadditivity of the Kodaira Dimension: Fibers of General Type," Advanced Studies in Pure Mathematics 10,
1987, ALgebraic Geometry, Senday, 1985, 361-398.






\bibitem{ko} S.J.  Kov\'acs, ``Smooth families over rational and elliptic curves," 
Jour.Alg.Geom. {\bf 5}  (1996), 369-385.


\bibitem{mi} L. Migliorini, ``A smooth family of minimal surfaces of general 
type over a curve of genus at most one is trivial,"
Jour.Alg.Geom. {\bf 4} (1995), 353-361.



\bibitem{mo2} S. Mori, ``Classification of higher dimensional varieties,"
in {\em Algebraic Geometry, Bowdoin 1985}, Proc. Symp. Pure Math., Vol. {\bf 46},
A.M.S. 1987.






\bibitem{pa} A.N. Parshin, ``Algebraic curves over function fields I,"
Isv. Akad. Nauk. SSSR, {\bf 32} (1968).




\bibitem{s-s} B. Shiffman, A.J. Sommese, {\em Vanishing Theorems
on complex manifolds}, Progress in Mathematics Vol. {\bf 56}, Birkh\"auser, Boston,
(1985)









\bibitem{za} O. Zariski, ``The theorem of Riemann-Roch for high multiples
of an effective divisor on an algebraic surface,"  Ann. of  Math., {\bf 76} (1962),
560-615.

\bibitem{qz} Q. Zhang, ``Global holomorphic one-forms on 
projective manifolds with ample canonical bundles," J. Algebraic Geom.
{\bf 6} (1997), no.4, 777-787.

\end{thebibliography}
\end{document}